\documentclass[a4paper, 10pt, conference]{ieeeconf}      

\IEEEoverridecommandlockouts                              
\usepackage{epsfig} 
\usepackage{psfrag}
\usepackage{amsmath} 
\usepackage{amssymb}  

\usepackage{url}

\newtheorem{thm}{Theorem}
\newtheorem{prp}{Proposition}
\newtheorem{cor}{Corollary}
\newtheorem{ass}{Assumption}
\newtheorem{lem}{Lemma}
\newtheorem{alg}{Algorithm}
\newtheorem{rem}{Remark}
\newtheorem{defin}{Definition}

\newenvironment{pf}{\smallbreak\noindent{\it Proof. }}{\hfill$\Box$\smallbreak}

\usepackage{enumerate}

\newcommand{\overbar}[1]{\mkern 1.5mu\overline{\mkern-1.5mu#1\mkern-1.5mu}\mkern 1.5mu}
\newcommand{\reals}{\mathbf{R}}
\newcommand{\hilbert}{\mathcal{H}}
\newcommand{\id}{{\rm{Id}}}
\DeclareMathOperator*{\argmin}{argmin}

\usepackage{tikz}

\usetikzlibrary{calc}
\usetikzlibrary{decorations.pathreplacing}

\newcommand{\marker}[2]
{
  \pgfgettransformentries{\myxscale}{\@tempa}{\@tempa}{\myyscale}{\@tempa}{\@tempa}
  \draw[thick] ($#1+0.08*(1/\myxscale,1/\myyscale)$)--($#1-0.08*(1/\myxscale,1/\myyscale)$);
  \draw[thick] ($#1+0.08*(-1/\myxscale,1/\myyscale)$)--($#1-0.08*(-1/\myxscale,1/\myyscale)$);
}

\newcommand{\radius}[5]
{   
  \draw[#5] (#1,#2)--({#1+#3*cos(#4)},{#2+#3*sin(#4)});
}

\newcommand{\secant}[3]
{   
  \draw ({#1*cos(#2)},{#1*sin(#2)})--({#1*cos(#3)},{#1*sin(#3)});
}

\newcommand{\secantcontr}[4]
{   
  \draw ({#1*cos(#2)},{#1*sin(#2)})--({#4*#1*cos(#3)},{#4*#1*sin(#3)});
}

\newcommand{\secantalphapoint}[4]
{ 
  \coordinate (P) at
  ({(1-#4)*#1*cos(#2)+#4*#1*cos(#3)},{(1-#4)*#1*sin(#2)+#4*#1*sin(#3)});
  \marker{(P)}{0.02};
}

\newcommand{\secantalphapointtext}[6]
{ 
  \node[#5] at
  ({(1-#4)*#1*cos(#2)+#4*#1*cos(#3)},{(1-#4)*#1*sin(#2)+#4*#1*sin(#3)}) {#6};
}

\newcommand{\secantalphapointcontr}[5]
{ 
  \coordinate (P) at ({#5*((1-#4)*#1*cos(#2)+#4*#1*cos(#3))},{#5*((1-#4)*#1*sin(#2)+#4*#1*sin(#3))});
  \marker{(P)}{0.02};
}

\newcommand{\secantalphapointtextcontr}[7]
{ 
  \node[#5] at
  ({#7*((1-#4)*#1*cos(#2)+#4*#1*cos(#3))},{#7*((1-#4)*#1*sin(#2)+#4*#1*sin(#3))}) {#6};
}

\newcommand{\arccenter}[6]
{
  \draw[#6] ([shift=(#4:#3)]#1,#2) arc (#4:#5:#3);
}

\title{\LARGE \bf
  Tight Linear Convergence Rate Bounds for\\Douglas-Rachford Splitting and ADMM}

\author{Pontus Giselsson$^\star$\thanks{$^\star$ Department of Automatic Control, Lund University.
   Email: {\tt{pontusg@control.lth.se}}}}

\begin{document}

\maketitle
\thispagestyle{empty}
\pagestyle{empty}

\begin{abstract}

Douglas-Rachford splitting and the alternating direction method of
multipliers (ADMM) can be used to solve convex optimization problems
that consist of a sum of two functions.
Convergence rate estimates for these algorithms have received much
attention lately. In particular, linear convergence rates have been
shown by several
authors under various assumptions. One such set of assumptions is
strong convexity and smoothness of one of the functions in the
minimization problem. The authors recently provided a linear convergence rate
bound for such problems. In this paper, we show that this rate bound
is tight for many algorithm parameter choices.

\end{abstract}


\section{Introduction}

Douglas-Rachford splitting is an optimization algorithm that can solve
general convex composite optimization problems. The algorithm has its roots in
the 1950's \cite{DouglasRachford,PeacemanRachford}. In the late
1970's, it was shown
\cite{LionsMercier1979} how to use the 
algorithm to solve montone operator inclusion problems and convex
composite optimization problems. The alternating direction
method of multipliers (ADMM) can also solve composite optimization
problems. It was first presented in
\cite{Glowinski1975,Gabay1976}. Soon thereafter,
it was shown \cite{Gabay83} that ADMM is equivalent to
Douglas-Rachford splitting applied to the dual problem. 

General sublinear convergence rate estimates for these methods have just
recently been presented in the literature, see
\cite{DR_one_over_k_2012,even_lin_conv_2013,conv_split_schemes_Davis_2014}.
Under various assumptions, also
linear convergence rates can be established. In the paper by Lions and
Mercier \cite{LionsMercier1979}, a linear convergence rate was provided for
Douglas-Rachford splitting under (the equivalence of) strong
convexity and smoothness assumptions. Until recently, further linear
convergence rate 
results have been scarce. The last couple of years, however, several
linear convergence rate results for both Douglas-Rachford splitting
and ADMM have been presented. These include \cite{linConvADMM,Davis_Yin_2014}, in which linear
convergence rates for ADMM are presented under various assumptions.
In \cite{lin_conv_DR_mult_block_2013}, linear convergence rates are
established for multiple splitting ADMM.
In \cite{Panos_acc_DR_2014}, it is shown that for a specific class of problems, the
Douglas-Rachford algorithm can be interpreted as a gradient method of
a function named the Douglas-Rachford envelope. By showing strong
convexity and smoothness properties of the Douglas-Rachford envelope
under similar assumptions on the underlying problem, a linear
convergence rate is established based on gradient algorithm theory. 
Very recently \cite{laurent_ADMM} appeared and showed linear
convergence of ADMM under smoothness and strong convexity assumptions
using the integral quadratic constraints (IQC) framework. The rate is
obtained by solving a series of a small semi-definite programs.
Common for all these linear convergence rate bounds are that they are not
tight for the class of problems under consideration, see
\cite[Section~IV.B]{gisBoydTAC2014metric_select}.

In \cite{arvind_ADMM}, linear convergence of ADMM is established under
more general assumptions than the above.
However, the assumptions are more difficult to verify for a
given problem. Tightness
is verified for a 2-dimensional example
in the Euclidean case.
In \cite{GhadimiADMM}, linear convergence for ADMM on strongly convex quadratic
optimization problem with inquality constraints is established. This
rate improves on the rates presented in
\cite{LionsMercier1979,linConvADMM,Davis_Yin_2014,lin_conv_DR_mult_block_2013,Panos_acc_DR_2014,laurent_ADMM}.
In \cite{gisBoyd2014CDCprecondADMM}, the authors generalize, using a completely
different machinery, the
results in \cite{GhadimiADMM} and in
\cite{gisBoydTAC2014metric_select} the results are
further generalized. More specifically, \cite{gisBoydTAC2014metric_select}
generalizes the results in \cite{GhadimiADMM} in the following three ways; (i) a wider class
of problems is considered, (ii) rates for both Douglas-Rachford
splitting and ADMM are provided, and (iii) the results in \cite{gisBoydTAC2014metric_select} hold
for general real Hilbert spaces as opposed to the Euclidean space
only in \cite{GhadimiADMM}. For the restricted class of problems
considered in \cite{GhadimiADMM}, the
convergence rate bounds in \cite{gisBoydTAC2014metric_select} and
\cite{GhadimiADMM} coincide.

The contribution of this paper is that we show tightness of the convergence
rate bounds presented in \cite{gisBoydTAC2014metric_select} for the class of problems
under consideration and for many algorithm parameters. This is done by
formulating examples, both for Douglas-Rachford splitting and ADMM,
for which the linear convergence rate bounds are satisfied with
equality. Similar lower convergence rate bounds have been presented in
\cite{laurent_ADMM}. The bounds in this paper cover wider classes of
problems and are less conservative.

\section{Notation}

We denote by $\reals$ the set of real numbers, $\reals^n$ the
set of real column-vectors of length $n$. Further
$\overbar{\reals}:=\reals\cup\{\infty\}$ denotes the extended
real line. Throughout this paper $\hilbert$ denotes a real separable Hilbert
space. Its inner product is denoted by $\langle\cdot,\cdot\rangle$, the
induced norm by $\|\cdot\|$, and the identity operator by $\id$. The
indicator function for a set $\mathcal{X}$ is denoted by $\iota_{\mathcal{X}}$. Finally, the
class of closed, proper, and convex
functions $f~:~\hilbert\to\overbar{\reals}$ is denoted by
$\Gamma_0(\hilbert)$.

\section{Preliminaries}

In this section we present, well known
concepts, results, operators, and algorithms that 
will be extensively used in the paper.
\begin{defin}[Orthonormal basis]
An \emph{orthonormal basis} $\{\phi_i\}_{i=1}^K$ for a (separable) Hilbert space
$\mathcal{H}$ is an orthogonal basis, i.e.
$\langle\phi_i,\phi_j\rangle=0$ if $i\neq j$, where each basis vector
has unit length, i.e. $\|\phi_i\|=1$.
\end{defin}
Hereon, $\phi_i$ will denote elemtens of an orthonormal basis.
\begin{rem}
The number of elements in the basis (the cardinality) $K$ is equal to
the dimension of the corresponding
Hilbert space, which might be $\infty$. Also, by
definition of a basis, each element $x\in\hilbert$ can be (uniquely)
decomposed as $x=\sum_{i=1}^K \langle
x,\phi_i\rangle\phi_i$, see \cite[Proposition~3.3.10]{Willem_fcn_analysis}.
\end{rem}

The reason why we consider separable Hilbert spaces is the following
proposition which can be found, e.g., in \cite[Proposition~3.3.12]{Willem_fcn_analysis}.
\begin{prp}
A Hilbert space is separable if and only if it has an orthonormal basis.
\end{prp}

We will also make extensive use of the following two propositions that
are proven, e.g., in \cite[Propsition~3.3.10]{Willem_fcn_analysis} and
\cite[Propsition~3.3.14]{Willem_fcn_analysis} respectively.
\begin{prp}[Parseval's identity]
In separable Hilbert spaces $\hilbert$, the squarred norm of each
element $x\in\hilbert$ satisfies
\begin{align*}
\|x\|^2 = \sum_{i=1}^K|\langle x,\phi_i\rangle|^2.
\end{align*}
\end{prp}
\begin{prp}[Riesz-Fischer]
In separable Hilbert spaces $\hilbert$, the sequence
$\sum_{i=1}^\infty a_i\phi_i$ converges if and only if
$\sum_{i=1}^\infty a_i^2<\infty$. Then
\begin{align*}
\left\|\sum_{i=1}^K a_i\phi_i\right\|^2=\sum_{i=1}^K a_i^2.
\end{align*}
\end{prp}

\begin{defin}[Strong convexity]
  \label{def:strConv}
  A function $f\in\Gamma_0(\hilbert)$ is $\sigma$-\emph{strongly
    convex} if 
  \begin{align*}
    f(x)\geq f(y)+\langle u,x-y\rangle+\tfrac{\sigma}{2}\|x-y\|^2
  \end{align*}
  holds for all $x,y\in\hilbert$ and all $u\in\partial f(y)$.
\label{def:strConv}
\end{defin}
\begin{defin}[Smoothness]
  A function $f\in\Gamma_0(\hilbert)$ is
  $\beta$-smooth if it is differentiable and
  \begin{align}
    f(x)\leq f(y)+\langle \nabla f(y),x-y\rangle+\tfrac{\beta}{2}\|x-y\|^2
    \label{eq:fSmooth_cvx}
  \end{align}
  holds for all $x,y\in\hilbert$.
\label{def:smoothness}
\end{defin}
\begin{defin}[Proximal operators]
  The \emph{proximal operator} of a function
  $f\in\Gamma_0(\hilbert)$ is defined as
  \begin{align*}
    {\rm{prox}}_{\gamma f}(y) :=
    \argmin_x\left\{f(x)+\tfrac{1}{2\gamma}\|x-y\|^2\right\}.
  \end{align*}
  \label{def:proxOp}
\end{defin}
\begin{defin}[Reflected proximal operators]
  The \emph{reflected proximal operator} to $f\in\Gamma_0(\hilbert)$ is defined as
  \begin{align*}
    R_{\gamma f} := 2{\rm{prox}}_{\gamma f}-\id.
  \end{align*}
  \label{def:reflRes}
\end{defin}
\begin{defin}[Fixed-point]
A point $y\in\hilbert$ is a \emph{fixed-point} to the (single-valued)
operator $A~:~\hilbert\to\hilbert$ if 
\begin{align*}
y = Ay.
\end{align*}
The set of fixed-points to $A$ is denoted by ${\rm{fix}}A$.
\end{defin}
\begin{alg}[Generalized Douglas-Rachford splitting]
The generalized Douglas-Rachford splitting algorithm is given by the iteration
\begin{align}
z^{k+1} = (1-\alpha)\id+\alpha R_{\gamma g}R_{\gamma f}z^k
\label{eq:DRsplitting}
\end{align}
where $\alpha\in(0,1)$ and $\gamma>0$ are algorithm parameters.
\label{alg:DRsplitting}
\end{alg}
\begin{rem}
In the general case, $\alpha$ is restricted to the interval $(0,1)$. Under the
assumptions used in this paper, a larger $\alpha$ can be used as well,
see \cite{gisBoydTAC2014metric_select}.
\end{rem}

\section{Linear convergence rates}

In this section, we state the linear convergence rate results for
Douglas-Rachford and ADMM in \cite{gisBoydTAC2014metric_select}.
The paper \cite{gisBoydTAC2014metric_select} considers 
optimization problems of the form
\begin{align}
\begin{tabular}{ll}
minimize & $f(x)+g(\mathcal{A}x)$
\end{tabular}
\label{eq:prob}
\end{align}
where $x\in\mathcal{H}$, and $f$, $g$, and $\mathcal{A}$ satisfy the
following assumptions:
\begin{ass}~
\begin{enumerate}[(i)]
\item The function $f\in\Gamma_0(\mathcal{H})$ is $\sigma$-strongly convex
  and $\beta$-smooth.
\item The function $g\in\Gamma_0(\mathcal{K})$.
\item $\mathcal{A}~:~\mathcal{H}\to\mathcal{K}$ is a surjective bounded linear operator.
\end{enumerate}
\label{ass:prob}
\end{ass}
Under the additional assumption that $\mathcal{A} = \id$ (which
implies that $\mathcal{K}=\mathcal{H}$),
Douglas-Rachford splitting can be applied to solve \eqref{eq:prob}. It
enjoys a linear convergence rate, as shown in
\cite[Theorem 1]{gisBoydTAC2014metric_select}. This result is restated
here for convenience.
\begin{thm}
  Suppose that Assumption~\ref{ass:prob} holds and that $\mathcal{A}=\id$. Then
  the generalized 
  Douglas Rachford algorithm (Algorithm~\ref{alg:DRsplitting}) converges linearly
  towards a fixed-point $\bar{z}\in{\rm{fix}}(R_{\gamma f}R_{\gamma g})$ with at
  least rate
  $|1-\alpha|+\alpha\max\left(\tfrac{\gamma\beta-1}{\gamma\beta+1},\tfrac{1-\gamma\sigma}{1+\gamma\sigma}\right)$, i.e. 
  \begin{align*}
    \|z^{k+1}-\bar{z}\|\leq\left(|1-\alpha|+\alpha\max\left(\tfrac{\gamma\beta-1}{\gamma\beta+1},\tfrac{1-\gamma\sigma}{1+\gamma\sigma}\right)\right)^k\|z^{0}-\bar{z}\|
  \end{align*}
for any $\gamma > 0$ and
$\alpha\in(0,\tfrac{2}{1+\max\left(\tfrac{1-\gamma\sigma}{1+\gamma\sigma},\tfrac{\gamma\beta-1}{1+\gamma\beta}\right)})$.
  \label{thm:DR_lin_conv_subdiff}
\end{thm}
\begin{rem}
The bound on the rate in Theorem~\ref{thm:DR_lin_conv_subdiff} can be optimized
with respect to the algorithm parameters $\alpha$ and $\gamma$. The
optimal parameters are given by $\alpha=1$ and $\gamma =
\tfrac{1}{\sqrt{\beta\sigma}}$ which yields rate bound factor
$\tfrac{\sqrt{\beta/\sigma}-1}{\sqrt{\beta/\sigma}+1}$, see \cite[Proposition 
16]{gisBoydTAC2014metric_select}. 
\label{rem:opt_param}
\end{rem}

In the case where $\mathcal{A}\neq\id$, problem \eqref{eq:prob} can
be solved by applying
Douglas-Rachford splitting on the dual problem:
\begin{align}
\begin{tabular}{ll}
minimize & $d(\mu)+g^*(\mu)$
\end{tabular}
\label{eq:dual_prob}
\end{align}
where $g^*\in\Gamma_0(\mathcal{K})$, and $d\in\Gamma_0(\mathcal{K})$ is
defined as
\begin{align*}
d := f^*\circ (-\mathcal{A}^*).
\end{align*}
If the dual problem \eqref{eq:dual_prob} satisfies
Assumption~\ref{ass:prob} (with $d$ instead of $f$ and $g^*$ instead
of $g$),
Douglas-Rachford splitting can be applied to solve
\eqref{eq:dual_prob}, and Theorem~\ref{thm:DR_lin_conv_subdiff}
would guarantee a linear convergence rate. Since $g\in\Gamma_0(\mathcal{K})$, we have
$g^*\in\Gamma_0(\mathcal{K})$ \cite[Theorem 12.2]{Rockafellar}, and we have $\mathcal{A}$ in
Assumption~\ref{ass:prob}(iii) equal to $\id$ in \eqref{eq:dual_prob}. The
remaining assumption needed to apply Theorem~\ref{thm:DR_lin_conv_subdiff} is that
$d\in\Gamma_0(\mathcal{K})$ is strongly convex and smooth. Indeed,
this is the case as shown in \cite[Proposition
18]{gisBoydTAC2014metric_select}. This result is restated here for
convenience of the reader.
\begin{prp}
  Suppose that Assumption~\ref{ass:prob} holds. Then
  $d\in\Gamma_0(\mathcal{K})$ is $\tfrac{\|\mathcal{A}^*\|^2}{\sigma}$-smooth and
  $\tfrac{\theta^2}{\beta}$-strongly convex, where $\theta>0$ always exists and satisfies
  $\|\mathcal{A}^*\mu\|\geq \theta\|\mu\|$ for all
  $\mu\in\mathcal{K}$. 
\label{prp:dual_fcn_prop_gen}
\end{prp}

It is well known \cite{Gabay83} that Douglas-Rachford splitting applied to the
dual problem \eqref{eq:dual_prob} is equivalent to ADMM applied to the
primal problem \eqref{eq:prob}. Therefore, the linear convergence rate
obtained by applying Douglas-Rachford splitting to the dual problem
\eqref{eq:dual_prob} directly translates to a linear convergence rate
for ADMM. This linear convergence rate bound is stated in \cite[Corollary
2]{gisBoydTAC2014metric_select}, and restated here for convenience.

\begin{prp}
  Suppose that Assumption~\ref{ass:prob}
  holds and that generalized Douglas-Rachford is applied to solve the
  dual problem \eqref{eq:dual_prob}. Then the Douglas-Rachford
  splitting algorithm converges linearly
  towards a fixed-point $\bar{z}\in{\rm{fix}}(R_{\gamma d}R_{\gamma g^*})$ with at
  least rate
  $|1-\alpha|+\alpha\max\left(\tfrac{\gamma\hat{\beta}-1}{\gamma\hat{\beta}+1},\tfrac{1-\gamma\hat{\sigma}}{1+\gamma\hat{\sigma}}\right)$, i.e. 
  \begin{align*}
    \|z^{k+1}-\bar{z}\|\leq\left(|1-\alpha|+\alpha\max\left(\tfrac{\gamma\hat{\beta}-1}{\gamma\hat{\beta}+1},\tfrac{1-\gamma\hat{\sigma}}{1+\gamma\hat{\sigma}}\right)\right)^k\|z^{0}-\bar{z}\|
  \end{align*}
for any $\gamma > 0$ and
$\alpha\in(0,\tfrac{2}{1+\max\left(\tfrac{1-\gamma\hat{\sigma}}{1+\gamma\hat{\sigma}},\tfrac{\gamma\hat{\beta}-1}{1+\gamma\hat{\beta}}\right)})$,
where $\hat{\beta} = \tfrac{\|\mathcal{A}^*\|^2}{\sigma}$ and
$\hat{\sigma} = \tfrac{\theta^2}{\beta}$.
\label{prp:ADMM_lin_conv_subdiff}
\end{prp}
\begin{rem}
The parameters that optimize
  the convergence rate bound are $\alpha=1$ and $\gamma =
  \tfrac{1}{\sqrt{\hat{\beta}\hat{\sigma}}}=\tfrac{\sqrt{\beta\sigma}}{\sqrt{\|\mathcal{A}^*\|^2\theta^2}}$ and the
  linear convergence rate bound factor is
  $\tfrac{\sqrt{\kappa}-1}{\sqrt{\kappa}+1}$, where $\kappa = \tfrac{\hat{\beta}}{\hat{\sigma}}=
  \tfrac{\|\mathcal{A}^*\|^2\beta}{\theta^2\sigma}$, see \cite[Corollary
2]{gisBoydTAC2014metric_select}.
\end{rem}








\section{Tightness of rate bounds}

In this section, we will state examples that show tightness of the linear convergence rate
bounds in Theorem~\ref{thm:DR_lin_conv_subdiff} and
Proposition~\ref{prp:ADMM_lin_conv_subdiff} for many choices of algorithm
parameters.

\subsection{Primal Douglas-Rachford splitting}

To establish that the convergence rate bound provided in
\cite[Theorem 1]{gisBoydTAC2014metric_select} and restated in
Theorem~\ref{thm:DR_lin_conv_subdiff}
is tight, we consider a problem of the form \eqref{eq:prob} with
\begin{align}
\label{eq:f_def} f(x) &=\sum_{i=1}^K\tfrac{\lambda_i}{2}\langle
x,\phi_i\rangle^2,\\
\label{eq:g_def} g(x) &=0,\\
\label{eq:A_def} \mathcal{A}&=\id.
\end{align}
Here $\{\phi_i\}_{i=1}^K$ is an orthonormal basis for $\hilbert$,
$K$ is the dimension of the space $\mathcal{H}$
(possibly infinite), and $\lambda_i$ is either $\sigma$
or $\beta$. We denote the set of indices $i$ with $\lambda_i=\sigma$
by $\mathcal{I}_{\sigma}$ and the set of indices $i$ with
$\lambda_i=\beta$ by $\mathcal{I}_{\beta}$. We require that
$\mathcal{I}_{\sigma}\neq\emptyset$, that
$\mathcal{I}_{\beta}\neq\emptyset$, and we get that
$\mathcal{I}_{\sigma}\cap\mathcal{I}_{\beta}=\emptyset$ and
$\mathcal{I}_{\sigma}\cup\mathcal{I}_{\beta} = \{1,\ldots,K\}$.

First, we show that $f$ in \eqref{eq:f_def} is defined (finite) for all $x\in\mathcal{H}$,
even if $\mathcal{H}$ is infinite 
dimensional. Obviously $f(x)\geq 0$ for
all $x\in\mathcal{H}$. We also have for arbitrary
$x\in\mathcal{H}$ that
\begin{align*}
f(x) = \sum_{i=1}^K\frac{\lambda_i}{2}\langle
x,\phi_i\rangle^2 \leq \frac{\beta}{2}\sum_{i=1}^K\langle
x,\phi_i\rangle^2=\frac{\beta}{2}\|x\|^2<\infty
\end{align*}
where the last equality follows from Parseval's identity. Therefore,
the optimization problem \eqref{eq:prob} with $f$, $g$, and $\mathcal{A}$ as in
\eqref{eq:f_def}, \eqref{eq:g_def}, and \eqref{eq:A_def} respectively
is well defined also on infinite dimensional spaces.

Next, we show that $f\in\Gamma_0(\hilbert)$ satisfies Assumption~\ref{ass:prob}(i),
i.e., that $f$ is $\beta$-smooth and $\sigma$-strongly convex.
\begin{prp}
The function $f$, as defined in \eqref{eq:f_def} with
$\lambda_i=\sigma$ for $i\in\mathcal{I}_{\sigma}$ and $\lambda_i=\beta$ for
$i\in\mathcal{I}_{\beta}$, is $\sigma$-strongly convex and $\beta$-smooth.
\end{prp}
\begin{pf}
Since $\mathcal{H}$ has a orthonormal basis, each element
$x\in\mathcal{H}$ may be decomposed as $x=\sum_{i=1}^K\langle
x,\phi_i\rangle\phi_i$. We let $a_i = \langle x,\phi_i\rangle$
and $b_i = \langle y,\phi_i\rangle$, to get arbitrary
$x=\sum_{i=1}^K  a_i\phi_i\in\mathcal{H}$ and $y=\sum_{i=1}^K
 b_i\phi_i\in\mathcal{H}$. Then
\begin{align*}
&\frac{\beta}{2}\|x-y\|^2 = \frac{\beta}{2}\|\sum_{i=1}^K a_i\phi_i-\sum_{i=1}^K b_i\phi_i\|^2\\
&=\sum_{i=1}^K \frac{\beta}{2}( a_i- b_i)^2\geq \sum_{i=1}^K\frac{\lambda_i}{2}( a_i- b_i)^2\\
&=\sum_{i=1}^K\lambda_i\left(\frac{1}{2} a_i^2 -\frac{1}{2} b_i^2 -\langle
 b_i\phi_i, a_i\phi_i- b_i\phi_i\rangle\right)\\
&=f(x)-f(y)-\sum_{i=1}^K\langle
\lambda_i\langle y,\phi_i\rangle \phi_i,( a_i- b_i)\phi_i\rangle\\
&=f(x)-f(y)-\langle
\sum_{i=1}^K\lambda_i\langle y,\phi_i\rangle \phi_i,\sum_{i=1}^K( a_i- b_i)\phi_i\rangle\\
&=f(x)-f(y)-\langle
\nabla f(y),x-y\rangle
\end{align*}
where the second equality follows from Riesz-Fischer, the first inequalty holds
since $\beta\geq\lambda_i$ for all $i=1,\ldots,K$, the third equality follows by
expanding the square and noting that
$ a_i b_i=\langle a_i\phi_i, b_i\phi_i\rangle$ and
$ b_i^2=\langle b_i\phi_i, b_i\phi_i\rangle$, the 
fourth equality follows by identifying the definition of $f$
in \eqref{eq:f_def} and using $ a_i=\langle x,\phi_i\rangle$ and
$ b_i=\langle y,\phi_i\rangle$, the
fifth equality holds since the added cross-terms vanish in the inner
product expression due to orthogonality of basis vectors $\phi_i$,
and the final equality holds by identifying $x=\sum_{i=1}^K
 a_i\phi_i$, $y=\sum_{i=1}^K 
 b_i\phi_i$, and the gradient of $f$ in \eqref{eq:f_def}:
\begin{align*}
\nabla f(x) = \sum_i \lambda_i\langle x,\phi_i\rangle\phi_i.
\end{align*}
This is the definition of $\beta$-smoothness in Definition~\ref{def:smoothness}.

An equivalent derivation using $\sigma$ instead of $\beta$
and a reversed inequality, shows that $f$ is also $\sigma$-strongly
convex. 
\end{pf}

To show that the provided example converges exactly with the rate
given in Theorem~\ref{thm:DR_lin_conv_subdiff}, we need expressions for the proximal operators and
reflected proximal operators of $f$ and $g$ in \eqref{eq:f_def} and
\eqref{eq:g_def} respectively.
\begin{prp}
The proximal operator of $f$ in \eqref{eq:f_def} is
\begin{align}
&{\rm{prox}}_{\gamma
  f}(y)=\sum_{i=1}^K\tfrac{1}{1+\gamma\lambda_i}\langle y,\phi_i\rangle\phi_i
\label{eq:f_prox}
\end{align}
and the reflected proximal operator is
\begin{align}
R_{\gamma f}(y) &=\sum_{i=1}^K\tfrac{1-\gamma\lambda_i}{1+\gamma\lambda_1}\langle y,\phi_i\rangle\phi_i.
\label{eq:f_refl_prox}
\end{align}
\label{prp:f_refl_prox}
\end{prp}
\begin{pf}
We decompose $x = \sum_{i=1}^K  a_i\phi_i$ where $ a_i=\langle
x,\phi_i\rangle$ and $y=\sum_{i=1}^K  b_i\phi_i$ where $ b_i=\langle
y,\phi_i\rangle$. Then, for general $\gamma>0$, the proximal operator of $f$ is given by:
\begin{align*}
&{\rm{prox}}_{\gamma f}(y)=\arg\min_x\left\{\gamma \left(\sum_{i=1}^K\tfrac{\lambda_i}{2}\langle
  \phi_i,x\rangle^2\right)+\tfrac{1}{2}\|x-y\|^2\right\}\\
&=\arg\min_{x=\sum_{i=i}^K a_i\phi_i}\left\{
    \left(\sum_{i=1}^K\tfrac{\gamma\lambda_i}{2} a_i^2\right)+\frac{1}{2}\left\|\sum_{i=1}^K(a_i-b_i)\phi_i\right\|^2\right\}\\ 
&=\arg\min_{x=\sum_{i=i}^K a_i\phi_i}\left\{
    \frac{1}{2}\sum_{i=1}^K\left(\gamma\lambda_i a_i^2+\left( a_i- b_i\right)^2\right)\right\}\\ 
&=\sum_{i=1}^K\arg\min_{ a_i}\tfrac{1}{2}\left\{
    \gamma\lambda_i a_i^2+( a_i- b_i)^2\right\}\phi_i\\ 
&=\sum_{i=1}^K\tfrac{1}{1+\gamma\lambda_i} b_i\phi_i=\sum_{i=1}^K\tfrac{1}{1+\gamma\lambda_i}\langle
y,\phi_i\rangle\phi_i.
\end{align*}
The reflected resolvent for general $\gamma>0$ is given by:
\begin{align*}
R_{\gamma f}(y) &= 2{\rm{prox}}_{\gamma f}(y)-y\\
&=2\sum_{i=1}^K\tfrac{1}{1+\gamma\lambda_i} b_i\phi_i-\sum_{i=1}^K b_i\phi_i\\
&=\sum_{i=1}^K\tfrac{1-\gamma\lambda_i}{1+\gamma\lambda_i} b_i\phi_i=\sum_{i=1}^K\tfrac{1-\gamma\lambda_i}{1+\gamma\lambda_i}\langle
y,\phi_i\rangle\phi_i.
\end{align*}
\end{pf}

The proximal and reflected proximal operators of $g\equiv 0$ are trivially
given by ${\rm{prox}}_{\gamma g} = R_{\gamma g} = \id$.

Next, these results are used to show a lower bound on the
convergence rate of Douglas-Rachford splitting for several choices of
algorithm parameters $\alpha$ and $\gamma$. First, we state two help
lemmas.

\begin{lem}
The function $\psi(x)=\tfrac{1-x}{1+x}$ is a decreasing function for
$x>-1$.
\label{lem:psi_decreasing}
\end{lem}
\begin{pf}
We have
\begin{align*}
&&(1-x)/(1+x)&< (1-y)(1+y) \\
&\Leftrightarrow& (1-x)(1+y)&< (1-y)(1+x) \\
&\Leftrightarrow& 2y&< 2x.
\end{align*}
\end{pf}
\begin{lem}
For $x> -1$, the function $\psi(x)=\tfrac{1-x}{1+x}$ satisfies $\phi(x)\leq
-\phi(y)$ if and only if $y\geq 1/x$.
\label{lem:psi_reciprocal}
\end{lem}
\begin{pf}
We have
\begin{align*}
&&\phi(x) = (1-x)/(1+x) &\leq (y-1)(1+y) = -\phi(y)\\
&\Leftrightarrow& (1-x)(1+y) &\leq (y-1)(1+x)\\
&\Leftrightarrow& 2 &\leq 2xy.
\end{align*}
\end{pf}

\begin{thm}
The generalized Douglas-Rachford splitting algorithm 
(Algorithm~\ref{alg:DRsplitting}) when applied to solve \eqref{eq:prob}
with $f$, $g$, and $\mathcal{A}$ in \eqref{eq:f_def}-\eqref{eq:A_def}
converges exactly with rate 
\begin{align}
|1-\alpha|+\alpha\max\left(\tfrac{1-\gamma\sigma}{1+\gamma\sigma},\tfrac{\gamma
    \beta-1}{1+\gamma\beta}\right)
\label{eq:DRrate_tight}
\end{align}
in the following cases: (i) $\alpha = 1$ and $\gamma\in(0,\infty)$, (ii) $\alpha\in(0,1]$ and
$\gamma\in(0, \tfrac{1}{\sqrt{\sigma \beta}}]$, (iii) $\alpha \in
[1,\tfrac{2}{1+\max\left(\tfrac{1-\gamma\sigma}{1+\gamma\sigma},\tfrac{\gamma
    \beta-1}{1-\gamma \beta}\right)})$ and $\gamma
  \in[\tfrac{1}{\sqrt{\sigma \beta}},\infty)$, (iv) $\alpha\in(0,\tfrac{2}{1+\max\left(\tfrac{1-\gamma\sigma}{1+\gamma\sigma},\tfrac{\gamma
    \beta-1}{1-\gamma \beta}\right)})$ and $\gamma = \tfrac{1}{\sqrt{\beta\sigma}}$.
\label{thm:tight_bound}
\end{thm}
\begin{pf}
For algorithm initial condition $z^0=\phi_i$ the Douglas-Rachford
algorithm evolves according to
\begin{align*}
z^{k} = \left(1-\alpha+\alpha\tfrac{1-\gamma\lambda_i}{1+\gamma\lambda_i}\right)^k\phi_i
\end{align*}
where $\lambda_i$ is either $\sigma$ or $\beta$ depending on if
$i\in\mathcal{I}_{\sigma}$ or $i\in\mathcal{I}_{\beta}$.
This follows immediately from Algorithm~\ref{alg:DRsplitting}, the
expression of $R_{\gamma f}$ in 
Proposition~\ref{prp:f_refl_prox}, and since $R_{\gamma g}=\id$.
Obviously, this converges with rate factor
\begin{align*}
\left|1-\alpha+\alpha\tfrac{1-\gamma\lambda_i}{1+\gamma\lambda_i}\right|.
\end{align*}
Below, we show for each of the four cases that this rate coincides with the
rate \eqref{eq:DRrate_tight}.
\subsection*{Case (i): $\alpha = 1$ {\rm{and}} $\gamma\in(0,\infty)$} The rate in this case when
$z^0 = \phi_i$, $i\in\mathcal{I}_\sigma$, is exactly
$\left|\tfrac{1-\gamma\sigma}{1+\gamma\sigma}\right|$.
The rate when $z^0 = \phi_i$, $i\in\mathcal{I}_\beta$, is exactly
$\left|\tfrac{1-\gamma \beta}{1+\gamma \beta}\right|$. 
A lower bound on the convergence of the algorithm when $\alpha = 1$ is therefore
\begin{align*}
\max\left(\left|\tfrac{1-\gamma\sigma}{1+\gamma\sigma}\right|,\left|\tfrac{1-\gamma
    \beta}{1+\gamma \beta}\right|\right)&=\max\left(\tfrac{1-\gamma\sigma}{1+\gamma\sigma},\tfrac{\gamma
    \beta-1}{1+\gamma \beta}\right)\\
&=|1-\alpha|+\alpha\max\left(\tfrac{1-\gamma\sigma}{1+\gamma\sigma},\tfrac{\gamma
    \beta-1}{1+\gamma \beta}\right).
\end{align*}
where the first equality is due to Lemma~\ref{lem:psi_decreasing}, and
the second holds since $\alpha=1$. This proves the first claim. 

\subsection*{Case (ii): $\alpha\in(0,1]$ {\rm{and}} $\gamma\in(0,\tfrac{1}{\sqrt{\sigma
  \beta}}]$} 
The rate when using initial condition $z^0 = \phi_i$,
$i\in\mathcal{I}_\sigma$, is
$r_{\sigma}:=1-\alpha+\alpha\tfrac{1-\gamma\sigma}{1+\gamma\sigma}$ (since
$(1-\alpha)\geq 0$ and
$\alpha\tfrac{1-\gamma\sigma}{1+\gamma\sigma}\geq 0$). For $z^0 = \phi_i$,
$i\in\mathcal{I}_\beta$, and $\gamma \leq \tfrac{1}{\beta}$, we get
\begin{align*}
|1-\alpha+\alpha\tfrac{1-\gamma \beta}{1+\gamma \beta}|&\leq
|1-\alpha|+|\alpha\tfrac{1-\gamma \beta}{1+\gamma \beta}|\\&=
1-\alpha+\alpha\tfrac{1-\gamma \beta}{1+\gamma \beta}\\
&\leq 1-\alpha+\alpha\tfrac{1-\gamma\sigma}{1+\gamma\sigma}=r_\sigma
\end{align*}
where the last inequality holds due to Lemma~\ref{lem:psi_decreasing}.
For $z^0 = \phi_i$,
$i\in\mathcal{I}_\beta$, and $\gamma\in[\tfrac{1}{\beta},\tfrac{1}{\sqrt{\sigma \beta}}]$, we get
\begin{align*}
|1-\alpha+\alpha\tfrac{1-\gamma \beta}{1+\gamma \beta}|&\leq
|1-\alpha|+|\alpha\tfrac{1-\gamma \beta}{1+\gamma \beta}|\\&=
1-\alpha+\alpha\tfrac{\gamma \beta-1}{1+\gamma \beta} \\
&\leq 1-\alpha+\alpha\tfrac{1-\gamma\sigma}{1+\gamma\sigma} = r_\sigma
\end{align*}
where the last inequality follows from Lemma~\ref{lem:psi_reciprocal}. Thus, a lower bound on
the rate for $\alpha \in(0,1]$ and $\gamma \in(0,\tfrac{1}{\sqrt{\sigma \beta}}]$ is
\begin{align*}
r_\sigma=1-\alpha+\alpha\tfrac{1-\gamma\sigma}{1+\gamma\sigma}=
|1-\alpha|+\alpha\max\left(\tfrac{1-\gamma\sigma}{1+\gamma\sigma},\tfrac{\gamma
    \beta-1}{1+\gamma \beta}\right).
\end{align*}
This proves the second claim. 

\subsection*{Case (iii): $\alpha\in[1,\tfrac{2}{1+\max\left(\tfrac{1-\gamma\sigma}{1+\gamma\sigma},\tfrac{\gamma
    \beta-1}{1-\gamma \beta}\right)})$ {\rm{and}} $\gamma\in[\tfrac{1}{\sqrt{\sigma
    \beta}},\infty)$}
The rate when using $z^0 = \phi_i$,
$i\in\mathcal{I}_\beta$, is
$r_\beta:=\alpha-1+\alpha\tfrac{\gamma \beta-1}{1+\gamma \beta}$ (since
$(1-\alpha)\leq 0$ and
$\alpha\tfrac{1-\gamma \beta}{1+\gamma \beta}\leq 0$). For $z^0 = \phi_i$,
$i\in\mathcal{I}_\sigma$, and $\gamma\in[\tfrac{1}{\sqrt{\sigma \beta}},\tfrac{1}{\sigma}]$ the rate is
\begin{align*}
|1-\alpha+\alpha\tfrac{1-\gamma\sigma}{1+\gamma\sigma}|&\leq
|1-\alpha|+|\alpha\tfrac{1-\gamma \sigma}{1+\gamma \sigma}|\\&=
\alpha-1+\alpha\tfrac{1-\gamma \sigma}{1+\gamma \sigma}\\
&\leq \alpha-1+\alpha\tfrac{\gamma \beta-1}{1+\gamma \beta}=r_\beta
\end{align*}
where the last inequality follows from Lemma~\ref{lem:psi_reciprocal}. 
For $z^0 = \phi_i$,
$i\in\mathcal{I}_\sigma$, and $\gamma \geq \tfrac{1}{\sigma}$, we get
\begin{align*}
|1-\alpha+\alpha\tfrac{1-\gamma \sigma}{1+\gamma \sigma}|&\leq
|1-\alpha|+|\alpha\tfrac{1-\gamma \sigma}{1+\gamma \sigma}|\\&=
\alpha-1+\alpha\tfrac{\gamma \sigma-1}{1+\gamma \sigma}\\
&\leq \alpha-1+\alpha\tfrac{\gamma \beta-1}{1+\gamma \beta}=r_\beta
\end{align*}
where the last inequality is due to Lemma~\ref{lem:psi_decreasing}. This implies that a
lower bound on the rate for $\alpha \in[1,\tfrac{2}{1+\max\left(\tfrac{1-\gamma\sigma}{1+\gamma\sigma},\tfrac{\gamma
    \beta-1}{1-\gamma \beta}\right)})$ and $\gamma \in [\tfrac{1}{\sqrt{\sigma \beta}},\infty)$ is
\begin{align*}
r_\beta=\alpha-1+\alpha\tfrac{\gamma \beta-1}{1+\gamma \beta}=
|1-\alpha|+\alpha\max\left(\tfrac{1-\gamma\sigma}{1+\gamma\sigma},\tfrac{\gamma
    \beta-1}{1+\gamma \beta}\right).
\end{align*}

\subsection*{Case (iv): $\alpha\in(0,\tfrac{2}{1+\max\left(\tfrac{1-\gamma\sigma}{1+\gamma\sigma},\tfrac{\gamma
    \beta-1}{1-\gamma\beta}\right)})$ {\rm{and}} $\gamma=\tfrac{1}{\sqrt{\sigma
    \beta}}$}
This Case follows directly from Cases (ii) and (iii).
\end{pf}

The convergence rate for the example given by
$f$ and $g$ in \eqref{eq:f_def} and \eqref{eq:g_def} respectively
coincides with the upper bound on the convergence rate in 
\cite[Theorem 1]{gisBoydTAC2014metric_select} (which
is restated in Theorem~\ref{thm:DR_lin_conv_subdiff}). The
bound in \cite[Theorem 1]{gisBoydTAC2014metric_select} is
therefore tight for the class of problems under consideration and for the
combination of algorithm parameters specified in Theorem~\ref{thm:tight_bound}.

\begin{rem}
The upper bound on the rate in \cite[Theorem
1]{gisBoydTAC2014metric_select}, relies on the triangle
inequality between $(1-\alpha)(z^k-\bar{z})$ and $\alpha (R_{\gamma g}R_{\gamma
  f}z^k-R_{\gamma g}R_{\gamma f}\bar{z})$. To get equality, we must
find $\alpha$, $\gamma$ and $z^k$ such
that $(1-\alpha)(z^k-\bar{z})$ and $\alpha (R_{\gamma g}R_{\gamma
  f}z^k-R_{\gamma g}R_{\gamma f}\bar{z})$ are parallel. For remaining
combinations of $\gamma$ and $\alpha$, these become anti-parallel, and
the rate bound is not met exactly. Note, however, that for optimal
choices of $\alpha$ and $\gamma$, the bound is tight.
\end{rem}

\subsection{Dual Douglas-Rachford splitting (ADMM)}

This section concerns tightness of the rate bounds when
Douglas-Rachford splitting is applied to the dual problem
\eqref{eq:dual_prob}, or 
equivalently, when ADMM applied to the primal problem \eqref{eq:prob}.
To show tightness in this case, we consider the following problem
\begin{align}
\label{eq:f_def_dual} f(x) &=\sum_{i=1}^K\tfrac{\lambda_i}{2}\langle
x,\phi_i\rangle^2\\
\label{eq:g_def_dual}  g(x) &=\iota_{x=0}(x)\\
\label{eq:A_def_dual} \mathcal{A}(x) &= \sum_{i=1}^K\nu_i\langle x,\phi_i\rangle\phi_i
\end{align}
where $\lambda_i=\sigma$ and $\nu_i=\theta>0$ if $i\in\mathcal{I}_{\sigma}$ and
$\lambda_i=\beta$ and $\nu_i=\zeta>\theta$ if
$i\in\mathcal{I}_{\beta}$, where $\mathcal{I}_{\sigma}$ and
$\mathcal{I}_{\beta}$ are the same as before.
That $\mathcal{A}$ is linear follows trivially. That it is
self-adjoint, bounded, and surjective is shown in the following
proposition.
\begin{prp}
The linear operator $\mathcal{A}$ defined in \eqref{eq:A_def_dual} is self-adjoint, i.e.
$\mathcal{A}=\mathcal{A}^*$, and for every $x\in\mathcal{H}$, we have
\begin{align}
\theta\|x\|\leq\|\mathcal{A}(x)\|\leq\zeta\|x\|.
\label{eq:A_bounds}
\end{align}
Further $\|\mathcal{A}\|=\|\mathcal{A}^*\|=\zeta$.
\end{prp}
\begin{pf}
We start by showing that $\mathcal{A}$ is self-adjoint. We have
\begin{align*}
\langle\mathcal{A}(x),\mu\rangle &=\left\langle \sum_{i=1}^K\nu_i\langle
x,\phi_i\rangle\phi_i,\sum_{i=1}^K\langle\mu,\phi_i\rangle\phi_i\right\rangle\\
&=\sum_{i=1}^K\langle \nu_i\langle
x,\phi_i\rangle\phi_i,\langle\mu,\phi_i\rangle\phi_i\rangle\\
&=\sum_{i=1}^K\langle \langle
x,\phi_i\rangle\phi_i,\nu_i\langle\mu,\phi_i\rangle\phi_i\rangle\\
&=\left\langle \sum_{i=1}^K\langle
x,\phi_i\rangle\phi_i,\sum_{i=1}^K\nu_i\langle\mu,\phi_i\rangle\phi_i\right\rangle\\
&=\langle x,\mathcal{A}(\nu)\rangle
\end{align*}
where moving of summations are due to orthogonality of $\phi_i$. Next
we show the first inequality in \eqref{eq:A_bounds}:
\begin{align*}
\|\mathcal{A}(x)\| &=
\left\|\theta\sum_{i\in\mathcal{I}_{\sigma}}\langle
  x,\phi_i\rangle\phi_i+\zeta\sum_{i\in\mathcal{I}_{\beta}}\langle
  x,\phi_i\rangle\phi_i\right\|\\
&\geq\theta\left\|\sum_{i=1}^K\langle
  x,\phi_i\rangle\phi_i\right\|=\theta\|x\|
\end{align*}
since $0<\theta\leq\zeta$. The second inequality in
\eqref{eq:A_bounds} is proven similarly.
Finally, we show $\|\mathcal{A}\|=\zeta$. We have already shown that
$\|\mathcal{A}(x)\|\leq\zeta\|x\|$ for all $x\in\mathcal{H}$, i.e.,
that $\|\mathcal{A}\|\leq\zeta$. By definition of the operator
norm, we also know that $\|\mathcal{A}\|\geq \|\mathcal{A}(x)\|$ for
all $x\in\mathcal{H}$ with $\|x\|\leq 1$. Choosing $x=\phi_j$ (which
satisfies $\|x\|=\|\phi_j\|=1$) for any
$j\in\mathcal{I}_{\beta}$ (i.e. $j$ with $\nu_j=\zeta$) gives
\begin{align*}
\|\mathcal{A}\|\geq \|\mathcal{A}(\phi_j)\| = \left\|\sum_{i=1}^K\nu_i\langle \phi_j,\phi_i\rangle\right\|=\|\nu_j\|=\zeta.
\end{align*}
Thus, $\|\mathcal{A}\|=\zeta$ and the proof is complete. 
\end{pf}

This result implies that the assumptions in \cite[Corollary
2]{gisBoydTAC2014metric_select} (and
Proposition~\ref{prp:ADMM_lin_conv_subdiff}) are met by 
$f$, $g$, and $\mathcal{A}$ in \eqref{eq:f_def_dual},
\eqref{eq:g_def_dual}, and \eqref{eq:A_def_dual} respectively. The
bound on the convergence rate from \cite[Corollary
2]{gisBoydTAC2014metric_select} (and restated in
Proposition~\ref{prp:ADMM_lin_conv_subdiff}) is therefore valid.
To show that this bound is tight for
the class of problems under consideration, we need the following
explicit characterization of $d$:
\begin{align*}
d(\mu) &:= f^*(-\mathcal{A}^*\mu) = f^*(-\mathcal{A}\mu)\\
&=\sup_{x}\left\{\langle -\mathcal{A}\mu,x\rangle-f(x)\right\}\\
&=-\inf_{x}\left\{f(x)+\langle\mathcal{A}\mu,x\rangle\right\}\\
&=-\inf_{x}\left\{\sum_{i=1}^K\tfrac{\lambda_i}{2}\langle
  x,\phi_i\rangle^2+\langle
  \sum_{i=1}^K\nu_i\langle\mu,\phi_i\rangle\phi_i,x\rangle\right\}\\
&=-\inf_{ a_i}\Bigg\{\sum_{i=1}^K\tfrac{\lambda_i}{2}\langle
  \sum_{i=1}^K a_i\phi_i,\phi_i\rangle^2\\
&\qquad\qquad\qquad\qquad+\langle
\sum_{i=1}^K\nu_i\langle\mu,\phi_i\rangle\phi_i,\sum_{i=1}^K a_i\phi_i\rangle\Bigg\}\\
&=-\sum_{i=1}^K\inf_{ a_i}\left\{\tfrac{\lambda_i}{2}\langle
   a_i\phi_i,\phi_i\rangle^2+\langle
\nu_i\langle\mu,\phi_i\rangle\phi_i, a_i\phi_i\rangle\right\}\\
&=-\sum_{i=1}^K\inf_{ a_i}\left\{\tfrac{\lambda_i}{2}
   a_i^2+\nu_i\langle\mu,\phi_i\rangle a_i\right\}\\
&=-\sum_{i=1}^K\left\{\frac{(\nu_i\langle\mu,\phi_i\rangle)^2}{2\lambda_i}
  -\frac{(\nu_i\langle\mu,\phi_i\rangle)^2}{\lambda_i}\right\}\\
&=\sum_{i=1}^K\frac{(\nu_i\langle\mu,\phi_i\rangle)^2}{2\lambda_i}=\sum_{i=1}^K\tfrac{\nu_i^2}{2\lambda_i}\langle\nu,\phi_i\rangle^2
\end{align*}
where the decomposition $x=\sum_{i=1}^K a_i\phi_i$ with $a_i=\langle
x,\phi_i\rangle$ is used, and the
optimal $a_i=-\nu_i\langle\mu,\phi_i\rangle/\lambda_i$.
The function $d$ has exactly the same structure as the function $f$
but with $\lambda_i$ in $f$ in \eqref{eq:f_def} replaced by $\nu_i^2/\lambda_i$
in $d$.
The function $g^*$ is, for all $\mu\in\mathcal{H}$, given by
\begin{align*}
g^*(\mu) &=
\sup_{x\in\mathcal{H}}\left\{\langle\mu,x\rangle-\iota_{x=0}(x)\right\}=\langle\mu,0\rangle=0.
\end{align*}
This implies that the dual problem
\eqref{eq:dual_prob} with $f$, $g$, and $\mathcal{A}$ specified in
\eqref{eq:f_def_dual}, \eqref{eq:g_def_dual}, and
\eqref{eq:A_def_dual} has exactly the same structure as the primal problem
\eqref{eq:prob} with $f$ and $g$ specified in \eqref{eq:f_def} and
\eqref{eq:g_def} respectively and with $\mathcal{A}=\id$ . The only things
that differ are the scalars that 
multiply the quadratic terms in the functions $f$ and $d$
respectively. Therefore, we can immediately state the following corollary
to Theorem~\ref{thm:tight_bound}.
\begin{cor}
Let $f$ be given by \eqref{eq:f_def_dual}, $g$ be given by \eqref{eq:g_def_dual}, and $\mathcal{A}$ be
given by \eqref{eq:A_def_dual}. Then the generalized Douglas-Rachford
algorithm applied to solve the dual problem \eqref{eq:dual_prob} (or
equivalently ADMM applied to solve \eqref{eq:prob})
converges as in Theorem~\ref{thm:tight_bound} with $\beta$ and
$\sigma$ in Theorem~\ref{thm:tight_bound} replaced by
$\hat{\beta}=\tfrac{\|\mathcal{A}\|^2}{\sigma}$ and
$\hat{\sigma}=\tfrac{\theta^2}{\beta}$ respectively.
\label{cor:ADMM_lin_conv}
\end{cor}

The exact rate provided in Corollary~\ref{cor:ADMM_lin_conv},
coincides with rate bound 
in \cite[Corollary
2]{gisBoydTAC2014metric_select} and
Proposition~\ref{prp:ADMM_lin_conv_subdiff}. Therefore, we conclude that the rate
bound in \cite[Corollary
2]{gisBoydTAC2014metric_select} for ADMM on the primal problem, or equivalently for
Douglas-Rachford splitting on the dual problem, is tight for the class
of problems under consideration for many algorithm parameter choices.
Especially, the bound is tight for the
optimal paramters $\alpha$ and $\gamma$, as in the primal Douglas-Rachford case.


\section{Conclusion}

Recent results in the literature have shown linear convergence of
Douglas-Rachford splitting and ADMM under various assumptions. In this
paper, we have shown that the linear convergence rate bounds presented
in \cite{gisBoydTAC2014metric_select} are indeed tight for the class
of problems under consideration.

\bibliographystyle{plain}
\bibliography{/local/home/pontusg/Research/MPC/papers/references/references}

\begin{appendices}

\end{appendices}

\end{document}